\documentclass{amsart}
  \usepackage{amscd,amssymb,epsfig}
  \usepackage{graphicx}
   \usepackage{epic,eepic}

                     \numberwithin{equation}{subsection}

                     \newtheorem{propo}{Proposition}[subsection]
                     
                     \newtheorem{theor}[propo]{Theorem}
                     \newtheorem{lemma}[propo]{Lemma}
                     \theoremstyle{definition}

                     \theoremstyle{remark}

                     \newcommand{\ZZ}{\mathbb{Z}}

\newcommand{\A}{\mathcal{A}}

 \newcommand{\PPP}{\mathcal{P}}

\newcommand{\naws}{nanowords}

\newcommand{\N}{\mathcal{N}}

                     \newcommand{\id}{\operatorname{id}}

\begin{document}
      \title{Knots and  words}
                     \author[Vladimir Turaev]{Vladimir Turaev}
                     \address{%
              IRMA, Universit\'e Louis  Pasteur - C.N.R.S., \newline
\indent  7 rue Ren\'e Descartes \newline
                     \indent F-67084 Strasbourg \newline
                     \indent France \newline
\indent e-mail: turaev@math.u-strasbg.fr }
                     \begin{abstract}  Knots and links are interpreted as homotopy classes of nanowords and nanophrases  in an alphabet consisting of 4  letters. Similar results hold for  curves on surfaces. We also discuss versions of  the Jones link polynomial and the link quandles for nanophrases.
                     \end{abstract}
                     \maketitle

  \section {Introduction}

  C. F. Gauss \cite{ga} introduced a method allowing to encode closed
  planar curves by  words of a certain type called now Gauss words. This method     extends to  planar knot
   diagrams and their
  isotopy. This gives a   description of the set of isotopy types of
  classical knots in terms of   words and their transformations.

To state our   results we need to generalize both knots and words.
An appropriate generalization of knots is provided by  long virtual
knots, see \cite{kau}, \cite{gpv}. We shall use an equivalent
formulation in terms of stable equivalence classes of (pointed
oriented) knot diagrams on surfaces, see \cite{kk}, \cite{cks}. The
set of these equivalence classes $\mathcal K$ contains the set
$\mathcal S$ of isotopy classes of oriented  knots in $S^3$. On the
combinatorial side, the Gauss words generalize to so-called
nanowords, see
   \cite{tu2}.  Our main result
    is a  bijection
  between $\mathcal K$ and the set of (appropriately defined) homotopy classes
  of nanowords in the alphabet consisting of 4 letters. Note that the image of $\mathcal  S \subset \mathcal  K$ under this bijection
  can be described  via the Rosenstiehl  theorem \cite{ro} giving necessary and sufficient conditions for the
curve corresponding to a given Gauss word to be planar.  Thus,
isotopy
  classification of   knots in $S^3$ is a special instance of  homotopy classification of words.
  Similarly, classical and virtual  links can be interpreted as nanophrases.
  This gives a broader perspective to knot theory. A number of    methods of
   knot theory  including the Kauffman bracket polynomial, the Jones polynomial, the knot quandle etc.  can be extended to the more general setting of words and phrases in arbitrary alphabets.

The theory of knots and links being very reach, it is natural to
study  related   simpler objects. One simplification of links in
cylinders over surfaces is obtained by projecting them to the
surfaces, i.e.,  by forgetting the over/under-crossing information
in
  link diagrams. The   stable equivalence classes of  curves on
surfaces were studied in \cite{kad}, \cite{tu1}. We describe these
classes   in terms of nanowords and nanophrases over an alphabet
consisting of only two letters.  We also suggest two different
simplifications of knot theory in terms of
  nanophrases over a 2-letter alphabet. The resulting  combinatorial
  objects are called pseudo-links and
quasi-links.  Every oriented link in a cylinder over an
oriented surface projects to a (multi-component) curve on the
surface, to a pseudo-link and to a quasi-link. Their distinctive
features  are contained in the following   facts: the Jones
polynomial of a link depends only on the underlying pseudo-link; the
fundamental group of the 2-fold branched covering of a classical
link depends only on the underlying quasi-link.

The plan of the paper is as follows. In Sect.\ 2 we recall
  nanowords  and their homotopy. In Sect.\ 3 we    identify the  stable equivalence classes
of pointed curves on surfaces with homotopy classes of
                     nanowords in a 2-letter  alphabet. In Sect.\
                     4  we identify
the stable equivalence classes of pointed knot diagrams on surfaces
with homotopy classes of  nanowords in a 4-letter  alphabet. In
Sect.\ 5 we show how to get rid of the base points. In Sect.\ 6 we
extend these results to links. In Sect.\ 7 we introduce pseudo-links
and    quasi-links. In Sect.\ 8 we   discuss   the bracket
polynomial and the Jones polynomial.  In Sect.\ 9 and 10 we discuss the
keis of nanophrases.

Conventions. Throughout the paper, all surfaces, curves, knots,
links, and knot and link diagrams are oriented, unless explicitly
stated to the contrary.

     \section{Nanowords and  homotopy}

\subsection{Words and nanowords}\label{word}   An {\it alphabet}  is a    set and
{\it letters} are its elements.     A {\it   word of length} $n\geq
1$  in an alphabet  $\A$ is  a mapping $w: \hat n\to \A  $ where
$\hat n=\{1,2,...,n\}$.    A word $w: \hat n\to \A  $
 is usually encoded by the  sequence of letters $w(1) w(2) \cdots  w(n)$. A word  $w: \hat n\to \A  $ is a {\it Gauss word}
  if each element of $\A$ is
the image of precisely two elements of  $\hat n$.

For a set $\alpha$, an {\it $\alpha$-alphabet} is a set $\A$ endowed
with a mapping $ \A\to \alpha$ called {\it projection}. The image of
  $A\in \A$ under this mapping is denoted $\vert A\vert$. A
{\it nanoword}
  over $\alpha$ is  a pair (an $\alpha$-alphabet $\A$,
a Gauss word in the alphabet $\A$). For example, any Gauss word $w$
in the alphabet $\alpha$ yields a nanoword $(\A=\alpha,w)$, see
\cite{tu2} for more examples and further details.  By definition,
there is a unique {\it empty nanoword} $\emptyset$ of length 0.

  An {\it
isomorphism} of $\alpha$-alphabets $\A_1$,  $\A_2$ is a bijection
$f:\A_1\to \A_2$ such that    $\vert A\vert=\vert f(A)\vert$ for all
$A\in \A_1$. Two nanowords    ($\A_1$, $w_1$) and  ($\A_2$, $w_2$)
over $\alpha$  are {\it   isomorphic} if    there is an isomorphism
of $\alpha$-alphabets    $f:\A_1\to \A_2$
   such that  $w_2  =f  w_1$.

\subsection{Homotopy of nanowords}\label{honawo} A {\it homotopy data}
 consists of a set $\alpha$ with involution $\tau:\alpha\to
\alpha$ and a set $S\subset \alpha^3=\alpha\times \alpha \times
\alpha$. The following three transformations of {\naws} over
$\alpha$ are called {\it $S$-homotopy moves} or simply {\it homotopy
moves}.

(1). The first move  applies to any nanoword   of the form    $(\A,   x A
A  y)$ where $A\in \A $ and $x,y$ are words  in the alphabet $\A'=\A-\{A\}$. It transforms $(\A,   x A
A  y)$  into
the nanoword
$(\A',  xy)$ where  the structure of an $\alpha$-alphabet in $\A' $ is obtained by restricting the  one  in  $\A$.  Note that
$ x  y$ is   a Gauss word in the alphabet $\A'$.

 The inverse move
$(\A',   xy)\mapsto (\A'\cup \{A\} ,   x A  A  y) $    adds a  new
letter $A $   with   arbitrary  $\vert A\vert \in \alpha$ and
replaces the  Gauss word $xy$ in  the $\alpha$-alphabet $\A'$ with
$xAAy$.

(2). The second move  applies to a  nanoword  of the form $(\A ,   xA  B  y BAz  )$   where $A,B\in  \A$ with
$\vert B\vert =\tau ( \vert A\vert)$ and $x,y,z$ are words in the alphabet  $\A'=\A-\{A,B\}$.  This  nanoword is transformed   into
$(\A' ,   xyz)$ where  the structure of an $\alpha$-alphabet in $\A' $ is obtained by restricting
the  one in  $\A$.

The inverse move  $(\A', xyz)\mapsto  (\A'\cup \{A,B\},   xA B  y
BAz)$      adds   two  new letters $A , B$   with   arbitrary $\vert
A\vert \in \alpha$ and       $\vert B\vert =\tau ( \vert A\vert)$
and    replaces the Gauss word $xyz$ in  the alphabet $\A'$ with
$xA B  y BAz$.

(3)   The third move  applies to a  nanoword  of the form $(\A, xAByACzBCt)$   where  $A,B,C \in \A$ are distinct letters
 such that
$(\vert A\vert ,\vert B\vert , \vert C\vert)\in S $ and $x,y,z,t$
are words in the alphabet $\A-\{A,B,C\} $. The
  move  transforms  $(\A, xAByACzBCt)$  into  $ (\A, xBAyCAzCBt)$. The inverse move
applies if   $(\vert A\vert ,\vert B\vert , \vert C\vert)\in S $ and
transforms  a nanoword $ (\A, xBAyCAzCBt)$ into $ (\A, xAByACzBCt)$.

  Two
nanowords over $\alpha$  are {\it $S$-homotopic} if they can be
obtained from each other by a finite sequence of    homotopy moves
(1) -- (3), the inverse moves, and     isomorphisms.  The relation
of $S$-homotopy is denoted   $\simeq_S $.    The set of $S$-homotopy
classes of nanowords over  $\alpha$ is denoted $\N(\alpha,S) $. (The
involution $\tau$ is omitted in this notation for shortness.)

Recall the following two lemmas from  \cite{tu2}, Sect.\ 3.2.

\begin{lemma}\label{2edcolbrad}   Let $A,B,C $ be distinct letters in  an $\alpha$-alphabet $ \A$ and let $x,y,z,t$ be words in the alphabet $\A-\{A,B,C\} $ such
that $xyzt$ is a Gauss word in this alphabet. Then

(i) $ (\A, xAByCAzBCt) \simeq_S (\A, xBAyACzCBt)$  for $(\vert
A\vert ,\tau (\vert B\vert), \vert C\vert )\in S  $;

(ii)
 $ (\A, xAByCAzCBt) \simeq_S  (\A, xBAyACzBCt)$  for $(\tau (\vert A\vert) ,\tau (\vert
B\vert), \vert C\vert  ) \in S $;

(iii) $ (\A, xAByACzCBt) \simeq_S  (\A, xBAyCAzBCt)$   for $ (\tau
(\vert A\vert),\vert B\vert ,\vert C\vert )\in
 S
 $.  \end{lemma}

A homotopy data $(\alpha, S)$ is   {\it
 admissible} if $S\cap (\alpha \times a \times a) \neq \emptyset$ for all $a\in
 \alpha$. For instance, if $S$ contains the diagonal
 $\{(a,a,a)\}_{a\in \alpha}$ of $\alpha^3$, then $(\alpha, S)$ is
 admissible.

\begin{lemma}\label{1edcolbrad}    Let $(\A, xAByABz)$ be a nanoword over $\alpha$ where  $A,B\in  \A$ with
$\vert B\vert =\tau ( \vert A\vert)$ and $x,y,z$ are words in the
alphabet $\A-\{A,B\}$.   If $(\alpha, S)$ is admissible, then $ (\A,
xAByABz) \simeq_S  ( \A -\{A,B\}, xyz)$.
\end{lemma}

  A {\it
morphism} $(\alpha,   S)\to (\alpha',   S')$  between two homotopy
data  is an equivariant   mapping $f:\alpha \to \alpha'$ such that
  $(f\times f \times f) (S)\subset S'$. Given   $f$, we can transform a nanoword $(\A, w )$ over
$\alpha$ into a nanoword $ (\A',w')$ over $\alpha'$  where $\A'=\A$
as sets, $w'=w$, and the projection $\A'\to \alpha'$ is the
composition of the projection $\A'=\A\to \alpha$    with $f$. This
transformation is compatible with homotopy and induces a monoid
homomorphism $\N(\alpha,S)\to \N(\alpha',S')$.

 \section{Curves versus  words}

\subsection{Curves}\label{virss} By a {\it curve}, we   mean
the image of a generic  immersion of an oriented circle into an
  oriented surface. The word \lq\lq generic" means that the
curve has only a finite set of self-intersections which are all
double and transversal.  A curve is {\it pointed} if it is endowed
with a base point (the origin) which is not a self-intersection. Two
pointed curves are {\it stably homeomorphic} if there is a
homeomorphism of their regular neighborhoods in the  ambient
surfaces mapping the first curve onto the second one and preserving
the origin of the curve and the orientations of the curve and the
surface. In particular,  attaching   a 1-handle to the ambient
surface away from a curve or removing such a handle does not change
the stable homeomorphism type of the curve.

Following \cite{kk}, \cite{cks}, we call two pointed curves   {\it
stably equivalent} if they can be related by a finite sequence of
the following transformations: (i) replacing the curve with a stably
homeomorphic one; (ii) homotopy of the curve in  its ambient surface
away from the origin. Note that such a homotopy may push a branch of
the curve across another branch or a double point  but not across
the origin of the curve.

Denote $\mathcal C$   the set of stable equivalence classes of
pointed curves. This set  is a monoid with multiplication defined by
connected sum at the origin. We are
  far from understanding the algebraic structure of
$\mathcal C$.  Several stable equivalence invariants of pointed
curves were introduced   in \cite{tu1}, \cite{sw} where   curves are
studied in terms of     virtual strings.

We  show now that
 the study of  $\mathcal C$ is an instance of   homotopy theory of words.

\subsection{Homotopy data $(\alpha_0,   S_0)$}\label{0triplerss}
Consider the homotopy data $(\alpha_0,   S_0)$ where $\alpha_0$ is
the set   $\{a,b\}$ with involution $\tau:\alpha_0\to \alpha_0$
permuting $a,b$ and
 $S_0=\{(a,a,a)\}_{a\in \alpha_0}$ is the diagonal. This homotopy data is admissible in the sense
  of Sect.\ \ref{honawo}.

   \begin{theor}\label{dcppmolbrad}  There is a canonical bijection  $\mathcal C= \N(\alpha_0,S_0)$. \end{theor}
    \begin{proof} We  associate with any
    pointed curve $f$ a nanoword $ w(f)$ over  $\alpha_0$. Let us
    label  the double points of $f$  by (distinct) letters
    $A_1,...,A_n$ where $n$ is the number of double points.
   Starting at the origin of $f$ and following along $f$ in the positive direction  we write
   down the labels of all double points until the   return to the
   origin. Since every double point is traversed twice, this gives
   a Gauss word $ w(f)$   in the alphabet $\A=\{A_1,...,A_n\}$. Let $ t_i^1$ (resp.\ $t_i^2$) be the tangent
   vector to $f$ at the crossing point labeled by $A_i$ appearing at the first (resp.\ second)
   passage through this crossing. Set $\vert A_i\vert =a $ if the
   pair $(t_i^1, t_i^2)$ is positively oriented and $\vert A_i\vert =b $ otherwise. This makes $\A$ into an $\alpha_0$-alphabet
   and makes $w$ into a nanoword over $\alpha_0$. This
   nanoword   is well defined up to
   isomorphism.

   We claim  that stably equivalent pointed curves give rise to
   $S_0$-homotopic  nanowords. Stable homeomorphisms of curves  preserve the nanoword up to isomorphism.
We need to show that a homotopy of a curve $f$ in its ambient
surface away from the origin does not change the $S_0$-homotopy
class of $w(f)$. Such a homotopy can be obtained by an
 ambient isotopy   and a finite sequence of    local
deformations   shown in Figure 1  and the inverse deformations.  It
is understood that all deformations in Figure 1 are effected away
from the origin of $f$.

\begin{figure}
\centerline{\includegraphics[width=9cm]{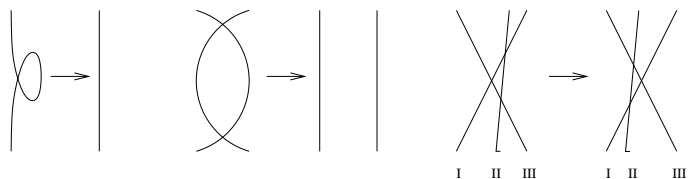}}
\caption{Three homotopy moves on curves}\label{figure1}
\end{figure}

 An ambient isotopy
does not change $w(f)$. A local deformation of the first type
changes $w(f)$ via the first homotopy move. Depending on the
orientations of the two branches of $f$, a local deformation of the
second type changes $w(f)$ via one of the  moves $  xA B  y BAz
 \mapsto  xyz $ or   $  xA B y ABz  \mapsto
 xyz $ where   $\vert B\vert =\tau ( \vert
A\vert)$. In both cases the $S_0$-homotopy class of $w(f)$ is
preserved. Consider the local deformation of the third type. It
suffices to consider the case where   all three branches
 are oriented   upwards. (The deformations  involving other orientations
of the branches can be obtained as compositions of this one with
ambient isotopy and local deformations of the second type.) There
are 6 cases to consider depending on the order in which one
traverses the three branches involved. Let $I$ (resp.\ $II,III$) be
the branch connecting the leftmost (resp.\ intermediate, rightmost)
bottom point to the rightmost (resp.\ intermediate, leftmost) top
point. If one traverses these branches in the order $I,II,III$
(resp.\ $III,II,I$), then $\vert A\vert =\vert B\vert =\vert C\vert
=a$ (resp.\ $\vert A\vert =\vert B\vert =\vert C\vert =b$) and the
deformation changes $w(f)$ via the third homotopy move (resp.\ its
inverse). If one traverses these branches in the order $II, III, I$
(resp.\ $I, III, II$), then the deformation changes $w(f)$ via the
second homotopy from Lemma \ref{2edcolbrad} (resp.\ its inverse).
Finally, if one traverses these three branches in the order $II,
I,III$ (resp.\ $III, I, II$), then the deformation changes $w(f)$
via the third homotopy from Lemma \ref{2edcolbrad} (resp.\ its
inverse).

We can conclude that the formula $f\mapsto w(f)$ defines a mapping
$W:\mathcal C\to  \N(\alpha_0, S_0)$. We claim that it is bijective.
It is easy to see that for every nanoword $w$ over $\alpha_0$ there
is a unique (up to stable homeomorphism) pointed curve $f$ such that
$w(f)=w$. Indeed, knowing $w$ we can uniquely recover an  oriented
regular
 neighborhood of such a curve in the ambient surface   (this well known construction    is described in
detail   in \cite{tu1}, Sect.\ 4.1  in terms of virtual strings.
Note that the notion of an open virtual string is equivalent to the
one of an isomorphism class of a  nanoword over $\alpha_0$.) This
implies that   $W$ is surjective.

  The injectivity of   $W$ follows
from the fact that if two nanowords are related by the $i$-th
homotopy move with $i=1,2,3$ then they can be represented by
  pointed curves related by the $i$-th deformation   in Figure
  1 (effected away from the origin). Here for $i=2$ the two branches  in Figure 1 are oriented in opposite
   directions  and for $i=3$
 all the branches   are oriented upwards and   traversed  in the order $I,II,III$.
    \end{proof}

\section{Knots versus  words}

\subsection{Knot diagrams}\label{2virss} By a {\it knot diagram}, we  mean
a (generic oriented) curve on an (oriented) surface such that at
each crossing point of the curve one of the two branches is
distinguished. The distinguished branch is  the {\it
over-crossing} and the second branch is   the {\it
under-crossing}. A knot diagram is {\it pointed} if it is endowed
with a base point (the origin) distinct from the crossing points.
Two pointed knot diagrams are {\it stably homeomorphic} if there is
a homeomorphism of their regular neighborhoods in the  ambient
surfaces  mapping the first diagram onto the second one and
preserving the origin,  the over/undercrossings, and the
orientations of the surface and the curve.

Following \cite{kk}, \cite{cks}, we call two pointed knot diagrams
{\it  stably equivalent} if they can be related by a finite sequence
of the following transformations: (i) replacing a knot diagram with
a stably homeomorphic one; (ii) the usual Reidemeister moves on a
knot diagram in its ambient surface away from the origin. The latter
moves  may push a branch of the diagram above or below a double
point or another branch but not across the origin. It should be
stressed that removing a closed subset from the ambient surface away
from a knot diagram or attaching a 1-handle away from the   diagram
does not change the stable equivalence type of the diagram.

Denote $\mathcal K$   the set of stable equivalence classes of
pointed knot diagrams.    The elements of $\mathcal K$ bijectively
correspond to long virtual knots in the sense of \cite{kau},
\cite{gpv}. Every  (oriented) knot   $K\subset S^3$ determines an
element of $\mathcal K$ obtained by presenting $K$ by a  diagram on
$S^2$ and picking an arbitrary base point. This yields a well
defined mapping from the set   of isotopy classes of   classical
knots  into $\mathcal K$. This mapping is essentially  injective,
see \cite{kau},
 \cite{gpv}.

Forgetting the over/under-crossing information we obtain  a natural
projection $\mathcal  K\to \mathcal  C$.  We now interpret $\mathcal
K$ in terms of words.

\subsection{Homotopy data $(\alpha_*,   S_*)$}\label{triplerss}
Consider the homotopy data $(\alpha_*,   S_*)$ where
$\alpha_*=\{a_+,a_-,b_+, b_-\}$ with involution $\tau:\alpha_*\to
\alpha_*$   defined by $\tau(a_{\pm})=b_{\mp},
\tau(b_{\pm})=a_{\mp}$ and $S_*\subset \alpha_*\times \alpha_*\times
\alpha_*$ consists of the following 12 triples:
$$ (a_{\pm},a_{\pm},a_{\pm}), (a_{\pm},a_{\pm},a_{\mp}),
(a_{\mp},a_{\pm},a_{\pm}),(b_{\pm},b_{\pm},b_{\pm}),
(b_{\pm},b_{\pm},b_{\mp}), (b_{\mp},b_{\pm},b_{\pm})  .$$ This
homotopy data is admissible in the sense
  of Sect.\ \ref{honawo}.

 Forgetting the signs, we obtain a projection
 $\alpha_*\to \alpha_0$. Applying it,  we can transform a
 nanoword over  $\alpha_*$ into a  nanoword over
 $\alpha_0$. This induces a monoid homomorphism
 $\N(\alpha_*,S_*)\to \N(\alpha_0,S_0)$.

   \begin{theor}\label{rt45pmolbrad}  There is a canonical bijection  $\mathcal K= \N(\alpha_*,S_*)$.
   Under this bijection, the monoid homomorphism
 $\N(\alpha_*,S_*)\to \N(\alpha_0,S_0)$ corresponds to the natural projection $\mathcal  K\to \mathcal  C$.
 \end{theor}
    \begin{proof} The proof reproduces the proof of Theorem \ref{dcppmolbrad} with a few changes.
    We begin by associating with any
    pointed knot diagram $F$ a nanoword $w=w(F)$ over  $\alpha_*$. As  usual,   each crossing of $F$ gives rise to
     a sign $\pm $.
      It is $+$ if the over-going branch crosses the under-going branch from left to right and $-$
      otherwise.
To define $w$,
    label  the double points of $F$  by (distinct) letters
    $A_1,...,A_n$ where $n$ is the number of double points.
   Starting at the origin of $F$ and following along $F$  we write
   down the labels of all double points until the   return to the
   origin. This gives
   a Gauss word $w$  in the alphabet $\A=\{A_1,...,A_n\}$. Let $ t_i^1$ (resp.\ $t_i^2$) be the tangent
   vector to $F$ at the crossing labeled $A_i$ appearing at the first (resp.\ second)
   passage through this crossing. Let $\varepsilon (i)=\pm$ be the sign of this crossing.
   Set $\vert A_i\vert =a_{\varepsilon (i)} $ if the
   pair $(t_i^1, t_i^2)$ is positively oriented and $\vert A_i\vert =b_{\varepsilon (i)} $ otherwise.
    This makes $\A$ into an $\alpha_*$-alphabet
   and makes $w=w(F)$ into a nanoword over $\alpha_*$. This
   nanoword   is well defined up to
   isomorphism.

   We need to verify that stably equivalent pointed knot diagrams give rise to
   $S_*$-homotopic  nanowords. Stable homeomorphisms preserve the nanoword up to isomorphism.   We need to show that
    the $S_*$-homotopy
class of $w(F)$ is preserved under the Reidemeister moves on $F$
away from the origin.
  The first Reidemeister move
changes $w(F)$ via the first homotopy move. Depending on the
orientations of the two branches of $F$ involved in the   second
Reidemeister move, the nanoword $w(F)$  changes   via one of the
moves $(\A , xA B  y BAz )\mapsto (\A-\{A,B\},xyz)$ or   $(\A , xA B
y ABz )\mapsto (\A-\{A,B\},xyz)$ where $A,B\in \A$ with $\vert
B\vert =\tau ( \vert A\vert)$. In both cases the $S_*$-homotopy
class of $w(F)$ is preserved. Consider the   third Reidemeister
move. It suffices to consider the case where   all three branches
 are oriented in the same direction, say upwards, and  the signs of all crossings are $+$. (This is the classical
  \lq\lq braid move" $\sigma_1 \sigma_2 \sigma_1\mapsto \sigma_2 \sigma_1 \sigma_2$; the
  moves involving other orientations
of the branches and/or other signs of crossings can be obtained as
compositions of this move with
  second Reidemeister moves.) There are 6 cases to
consider depending on the order in which one traverses the three
branches involved. Let $I$ (resp.\ $II$, $III$) be the branch
connecting the leftmost (resp.\ intermediate, rightmost) bottom
point to the rightmost (resp.\ intermediate, leftmost) top point. If
one traverses these branches in the order $I,II,III$ (resp.\
$III,II,I$), then $\vert A\vert =\vert B\vert =\vert C\vert =a_+$
(resp.\ $\vert A\vert =\vert B\vert =\vert C\vert =b_+$) and the
deformation changes $w(F)$ via the third homotopy move (resp.\ its
inverse) where we use that $(a_+,a_+, a_+)\in S_*$ (resp.\ that
$(b_+,b_+, b_+)\in S_*$). If one traverses these three branches in the
order $II, III, I   $ (resp.\ $I, III, II$), then the deformation
changes $w(F)$ via the second  homotopy from Lemma \ref{2edcolbrad}
(resp.\ its inverse) where we use that $(a_-, a_-,  a_+)\in S_*$
(resp.\ that $(b_-, b_-, b_+)\in S_*$). Finally, if one traverses
these branches in the order $(II,I,III)$ (resp.\ $III,I, II$), then
the deformation changes $w(F)$ via the third homotopy from Lemma
\ref{2edcolbrad} (resp.\ its inverse) where we use that $(b_+, b_-,
b_-)\in S_*$ (resp.\ that $(a_+,a_-, a_-)\in S_*$).

Thus the formula $F\mapsto w(F)$ defines a mapping $W:\mathcal K\to
\N(\alpha_*, S_*)$. We claim that $W$  is a bijection.   As in the
case of curves, for every nanoword $w$ over $\alpha_*$ there is a
unique (up to stable homeomorphism) pointed knot diagram $F$ such
that $w(F)=w$. Indeed it suffices to realize the underlying nanoword
over $\alpha_0$ by a   pointed curve and then to choose the
over/under-crossings to ensure the right signs at all crossings.
  This implies that   $W$ is surjective.

 To prove the injectivity
of   $W$
  it suffices to observe that if two nanowords over $\alpha_*$ are related
by the $i$-th homotopy move with $i=1,2,3$  then they can be
represented by
  pointed knot diagrams related by the $i$-th Reidemeister move (effected away from the origin).
  The cases $i=1,  2$ are straightforward. For $i=3$ the 12 elements of the set $S_*$ lead to   all 12    possible
  choices of over/under-crossings in the third move on Figure 1   leading to
  admissible
  Reidemeister
  moves (for this argument we  can assume  that
  the branches $I,II,III$   are oriented  upwards and traversed either  in the order $I,II,III$ or   $III,II,I$).

  The last claim of the theorem follows from the definitions.
    \end{proof}

\subsection{Examples}\label{45exallnseglss} The nanoword  $ABCABC$
with $\vert A\vert=\vert C\vert= a_+, \vert B\vert=b_+$ represents a  pointed 
trefoil, see Figure 2 where the thick point is the origin of the diagram. The nanoword  $ABCADCBD$ with $\vert A\vert=\vert D\vert=
b_+, \vert B\vert=b_-,\vert C\vert= a_-$ represents a  pointed   figure eight
knot, see Figure 2.

\begin{figure}
\centerline{\includegraphics[width=7cm]{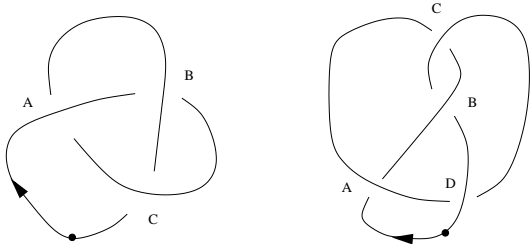}}
\caption{Trefoil and figure eight knot}\label{figure2}
\end{figure}

  \section{Eliminating the origin}\label{edro}

\subsection{Shifts}\label{funcyclomot} Fix an involution $\nu$ in a set $ \alpha$     called the  {\it shift involution}.
The {\it   $\nu$-shift} of a nanoword $(\A, w:\hat n \to \A)$ over $ \alpha$ is the
nanoword $(\A', w':\hat n \to \A')$
 obtained by moving the first letter $A=w(1) $ of $w$ to the end and applying
 $\nu$ to $\vert A\vert \in \alpha$. More precisely,     $\A' =(\A-\{A\}) \cup \{A_\nu\}$
  where $A_\nu$ is a \lq\lq new" letter not belonging to $\A$. The   projection  $ \A'\to \alpha$
   extends the   given  projection    $\A-\{A\}\to \alpha $ by $\vert A_\nu\vert=  \nu(\vert  A\vert)$.
    The word $w'$ in the alphabet $\A'$ is defined by $w'= x A_\nu y A_\nu$ for $w= AxAy
    $.

 Given a homotopy data $(\alpha,   S)$ and a shift involution $\nu$ in $\alpha$, we can quotient the set of nanowords over
 $\alpha$
by the equivalence relation generated by $S$-homotopy and
$\nu$-shifts. The resulting set is denoted $ {\N}(\alpha,S,\nu)$.
There is a natural projection $\N(\alpha,S) \to {\N}(\alpha,S,\nu)$
but there is no natural  multiplication in $ {\N}(\alpha,S,\nu)$.

\subsection{Non-pointed knots}\label{reffrss} Stable equivalence can be   defined for   (non-pointed)
knot diagrams through repeating the definition in the pointed case
but omitting all references to   base points. Denote $\hat{\mathcal
K}$ the set of stable equivalence classes of knot diagrams.  Each
knot in the cylinder over a  surface represents an element in
$\hat{\mathcal K}$ depending only on the isotopy type of the knot.

  As we know,
a pointed knot diagram gives rise to a nanoword in the alphabet
$\alpha_*=\{a_+,a_-,b_+, b_-\}$. This nanoword is preserved when the
origin is pushed along the generic part of the diagram. When the
origin jumps over a double point, the nanoword is modified by the
$\nu$-shift where  $\nu:\alpha_*\to \alpha_*$ is the involution
sending $a_\pm$ to $b_\pm$.   Theorem \ref{rt45pmolbrad} implies
that $\hat{\mathcal K}= {\N}(\alpha_*, S_*, \nu)$.

Non-oriented  knots can be treated similarly, we  do it  in the next
section in a more general setting of links.

 \section{Nanophrases and links}\label{naink}

\subsection{Nanophrases}\label{nanaopyclomot}   A {\it nanophrase}
 of length $k\geq 0$ over a set $\alpha$  is  a tuple consisting of  an $\alpha$-alphabet $\A$
  and a sequence of $k$ words $w_1,...,w_k$ in the alphabet $\A$ such
  that their concatenation $w_1w_2\cdots w_k$ is a Gauss word in
  this alphabet.  We   denote this nanophrase by  $(\A,
  (w_1\vert w_2 \vert \cdots \vert w_k))$ or shorter by $(w_1\vert w_2 \vert \cdots \vert
  w_k)$. Note that some of the words $w_1,...,w_k$ may be empty.

By definition, there is a unique {\it empty
nanophrase} of length 0 (the corresponding $\alpha$-alphabet $\A$ is void).

Any nanoword $w$ over $\alpha$ yields a nanophrase $ (w)$ of length
$1$. In the sequel  we make no difference between nanowords and
nanophrases of length 1. 

Isomorphism of two nanophrases is   an isomorphism of $\alpha$-alphabets transforming the first
sequence of words into the second one. Given a homotopy data $
\alpha, \tau, S$, we define homotopy moves on nanophrases as in
Sect.\ \ref{honawo} with the only difference that the 2-letter
sub-words $AA, AB, BA , AC, BC$ etc.\ modified by these moves    may
belong to different words of the phrase. Isomorphisms and homotopy
moves generate an equivalence relation $\simeq_S$ of $S$-homotopy on
the class of nanophrases over $\alpha$. Examples:
$$ (AB\vert AC\vert BC )\simeq_S (BA\vert CA\vert CB ),\,\,\,\, (AB\vert ADDC BC )\simeq_S (BA\vert CA  CB )$$ provided
$(\vert A\vert ,\vert B\vert ,\vert C\vert )\in S$. The length of a
nanophrase is preserved under   $S$-homotopy.

Lemmas \ref{2edcolbrad}    and \ref{1edcolbrad} extend to
nanophrases with the only change that the 2-letter sub-words   $AB,
BA, CA$ etc.\  may belong to different  words of the phrase.

\subsection{Operations on nanophrases}\label{opernanaopyclomot}
 Fix a homotopy data $(\alpha, \tau, S)$ and
a shift involution $\nu$ in $\alpha$. We define  $\nu$-shifts,
$\nu$-inversions, and $\nu$-permutations of words in a nanophrase
$P=(\A, (w_1\vert w_2 \vert \cdots \vert
  w_k))$ over
$\alpha$.

We can {\it $\nu$-shift} the $i$-th word $ w_i$ in $P $ through
moving the first letter, say  $A $, of $w_i$ to the end of $ w_i$
keeping $\vert A\vert \in \alpha$ if $A$ appears in $w_i$ only once
and applying $\nu$ to $\vert A\vert$ if $A$ appears in $w_i$ twice.
All other words in   $P$ are  preserved.

To define inversions, we need more notation. For a word $w$ in $\A$,
denote by $\A_w$ the same alphabet $\A$ with new projection $\vert
...\vert_w$ to $\alpha$ defined as follows: for $A \in \A$ set
$\vert A\vert_w=\tau (\vert A\vert)$ if $A  $ occurs in $w$ once,
$\vert A\vert_w=\nu (\vert A\vert)$ if $A $ occurs in $w$ twice, and
$\vert A\vert_w=
 \vert A\vert$ otherwise. The {\it $\nu$-inversion} of the
$i$-th word in   $ P$ replaces $w_i$ with the opposite  word $(w_i)^-$
obtained by reading  $w_i $ from right to left and replaces the
$\alpha$-alphabet $\A$ with $\A_{w_i}$. All other words in   $P$ are
preserved.

The words in  $P$ can  be permuted in an arbitrary way, producing
thus new nanophrases over $\alpha$. We   will need
   more sophisticated    permutations of words depending on    $\nu $.
We begin with notation. For two words $u,v$ in the  alphabet $\A$,
consider the mapping $\A\to \alpha$ sending $A\in \A $ to $\nu(\vert
A\vert) \in \alpha$ if $A$ appears both in $u$ and   $v$ and sending
$A $ to $ \vert A\vert $ otherwise. This mapping makes the set $\A $
into an $\alpha$-alphabet denoted $\A_{u \cap v}$. For
$i=1,...,k-1$, the {\it $\nu$-permutation} of the $i$-th and
$(i+1)$-st words transforms $P=(\A, (w_1\vert w_2 \vert \cdots \vert
  w_k))$ into the nanophrase $$(\A_{w_i\cap w_{i+1}}, (w_1\vert w_2 \vert \cdots
  \vert w_{i-1}\vert w_{i+1}\vert w_i\vert w_{i+2} \vert \cdots
  w_k)).$$ This operation
   is
  involutive. The $\nu$-permutations define  an action of the
  symmetric group  $S_k$ on the set of nanophrases   of
  length $k$.

Denote $ {\PPP} (\alpha,S,\nu)$ the set of nanophrases over
 $\alpha$ quotiented
by the equivalence relation generated by $S$-homotopy,
$\nu$-permutations and $\nu$-shifts on words. Denote $ {\PPP}_u
(\alpha,S,\nu)$ the set of nanophrases over
 $\alpha$ quotiented
by the equivalence relation generated by the same operations and the
$\nu$-inversions.

\subsection{Link diagrams}\label{linklomot}  {\it Link
diagrams} on   (oriented) surfaces are defined in the same way as
knot diagrams with the difference that they  may be formed by
several (transversal generic oriented closed) curves rather than
only one curve. These curves are   {\it components} of the diagram.
A link diagram is {\it pointed} if each    component  is endowed
with a base point (the origin) distinct from the crossing points of
the diagram. A link diagram is {\it ordered} if   its components are
numerated by $1,2,...,k$ where $k$ is the number of the components.
Two ordered pointed link diagrams are {\it stably homeomorphic} if
there is an orientation preserving homeomorphism of their regular
neighborhoods in the  ambient surfaces mapping the first diagram
onto the second one and preserving   the over/undercrossings  and
the order, the origins and the orientations of   the components.

The {\it  stable equivalence} of ordered pointed link diagrams is
generated by the same transformations as in the case of knots. These
transformations should  preserve the order and the origins of the
components; the Reidemeister moves are allowed only  away from the
origins.

Denote $\mathcal L$   the set of stable equivalence classes of
ordered pointed link diagrams. Recall the homotopy data $ \alpha_*,
\tau,  S_*$ defined in Sect.\ \ref{triplerss}.

   \begin{theor}\label{rtlinkslbrad}  There is a canonical bijection  $\mathcal L= \PPP(\alpha_*,S_*)$.
 \end{theor}

 The proof of this theorem is analogous to the proof of Theorem
 \ref{rt45pmolbrad}. To write down the nanophrase associated with an
 ordered pointed link diagram one   goes along the first component
 starting at its origin, then along the second component, etc.

  Forgetting  the order and
the origins of link components, we obtain a notion of stable
equivalence for (non-ordered non-pointed) link diagrams. Denote
$\hat {\mathcal L}$ the set of equivalence classes of link diagrams.
As in the case of knots, each   link   in the cylinder over a
  surface represents an element in $\hat{\mathcal L}$
depending only on the isotopy class of this link. Theorem
\ref{rtlinkslbrad} implies that $\hat {\mathcal
L}={\PPP}(\alpha_*,S_*,\nu)$  where $\nu:\alpha_*\to \alpha_*$ is
the involution sending $a_\pm$ to $b_\pm$.

Additionally forgetting   link orientations, we   obtain a notion of
stable equivalence for  unoriented link diagrams (on oriented
surfaces). Denote $\hat {\mathcal L}_u$ the set of equivalence
classes of unoriented link diagrams. Theorem \ref{rtlinkslbrad}
implies that $\hat {\mathcal L}_u={\PPP}_u (\alpha_*,S_*, \nu)$.

\subsection{Remarks} 1. Theorem  \ref{rtlinkslbrad} can be
extended to framed links.  Consider the involution $ \nu \tau =\tau
\nu:\alpha_*\to \alpha_*$ sending $a_\pm, b_\pm$ to $a_\mp, b_\mp$,
respectively. A {\it framed  homotopy} $\sim$ of nanophrases over
$\alpha_*$ is defined as the $S_*$-homotopy with the   first
homotopy move   replaced by the following \lq\lq framed homotopy
move" on a nanoword in a nanophrase:  $ x A A yBBz\mapsto xyz$
provided $\vert A\vert =\nu \tau (\vert B\vert)$. Framed homotopic
nanophrases are $S_*$-homotopic; the converse is in general not
true. Lemmas \ref{2edcolbrad} and \ref{1edcolbrad} for
$\alpha=\alpha_*$ extend to this setting by replacing  $\simeq_S$
with $\sim$. The proof of Lemma  \ref{2edcolbrad}   in \cite{tu2}
does not use the first homotopy move.  The proof of Lemma
\ref{1edcolbrad}     in \cite{tu2} uses the first homotopy move but
can be easily modified to use the framed move instead.

Let us show that a nanoword   $xAABBy$ is framed homotopic to   $xy$
provided $\vert A\vert =\tau (\vert B\vert)$. Pick letters $C,D$ not
appearing in $x,y$ with $\vert C\vert =\vert B\vert, \vert D\vert
=\vert A\vert$. Then
$$xAABBy  \sim   xACDABBCDy \sim xCADBACBDy\sim xDBBDy\sim
xy.$$ Here we   insert  $ CD \cdots CD $,  apply homotopy (iii) of
Lemma \ref{2edcolbrad},   delete  $ CA \cdots AC $, and finally
delete  $DBBD$. A similar argument shows that   a transformation $
xAABy\mapsto xBAAy$ preserves the
  framed homotopy class. These   observations easily imply that framed
  homotopy classes of nanophrases over $\alpha_*$   bijectively
  correspond to   stable equivalence classes of framed
ordered pointed link diagrams.

2.      The
results of Sect.\ \ref{edro} and \ref{naink} have an obvious version
for systems of transversal curves  on   surfaces; it suffices to replace   $\alpha_*$   by    $\alpha_0=\{a,b\}$.

\section{Pseudo-links and quasi-links}\label{erdt}

\subsection{Pseudo-links}\label{quasinanaopyclomot}  Set
$\alpha_1=\{1, -1\}$ with involution $\tau $ permuting $1$ and $-1$
and let $S_1\subset \alpha_1\times \alpha_1\times \alpha_1$ consist
of the following 6 triples:
$$ ( 1,1,1), (1,1,-1),(-1,1,1), (-1,-1,-1), (-1, -1, 1), (1, -1,-1)
 .$$ This
homotopy data is admissible in the sense
  of Sect.\ \ref{honawo}. As a shift involution in $\alpha_1$, we take the identity mapping $\id:\alpha_1\to \alpha_1$.
  The corresponding   permutations and shifts of words in
  nanophrases over $\alpha_1$ are the   ordinary permutations and cyclic
  shifts of words (involving no modification of the underlying
  $\alpha_1$-alphabets).
Nanophrases  over $\alpha_1$ considered up to permutations and
cyclic shifts of words are called
  {\it
  pseudo-links}.

This terminology is justified by  the following connections to knot
theory. Consider the projection $ \alpha_*\to \alpha_1 $ sending
$a_+, b_+$ to $1$ and $a_-,
  b_-$ to $-1$. This projection transforms $S_*\subset (\alpha_*)^3$ into $S_1$. It commutes with   $\tau$
 and  with
  the shift involutions in $  \alpha_*,  \alpha_1 $ (the
  shift involution   in $\alpha_* $ is defined by
  $\nu(a_\pm)=b_\pm$).
Applying the projection $  \alpha_*\to \alpha_1 $, we can transform
any nanophrase over $\alpha_*$
  into a nanophrase over $\alpha_1$. Clearly,    $S_*$-homotopic nanophrases over
  $\alpha_*$ yield $S_1$-homotopic nanophrases over
  $\alpha_1$. This induces a  mapping  ${\PPP}  (\alpha_*,S_* )\to {\PPP}  (\alpha_1,S_1 )$. Quotienting by permutations
  and shifts of words we obtain a  mapping from
  ${\PPP}  (\alpha_*,S_*,
  \nu)$ to ${\PPP}  (\alpha_1,S_1, \id)$. Further quotienting by
  inversions of words we obtain a  mapping
  ${\PPP}_u  (\alpha_*,S_*,
  \nu)\to {\PPP}_u  (\alpha_1,S_1, \id)$. All these mappings are
  surjective.

By  Sect.\ \ref{naink},  a   link
  diagram  $D$ on a     surface  yields a nanophrase over $\alpha_*$. Projecting to $\alpha_1$, we obtain
  a nanophrase  $p(D)$ over $\alpha_1$. If $D$ is pointed and ordered, then $p(D)$ is well-defined, otherwise $p(D)$ is
   defined only up
   to
  permutations and shifts of words. The class of $p(D)$ in ${\PPP}(\alpha_1,S_1, \id)$ is an invariant of stable
   equivalence of $D$. We call $p(D)$   the {\it underlying pseudo-link} of
  $D$. Further projecting to
   ${\PPP}_u (\alpha_1,S_1, \id)$ we obtain an invariant
   independent of the orientation of $D$.

   In the next section we   explain
  that the   Jones
  polynomial of a link  depends only on the underlying pseudo-link.
  This shows that   pseudo-links are highly non-trivial objects   retaining important features of   links.

\subsection{Quasi-links}\label{quasinanaopyclomot}  Set
$\alpha_2=\{c, d\}$ with the identity involution $\tau
=\id:\alpha_2\to \alpha_2$   and let $S_2\subset \alpha_2\times
\alpha_2\times \alpha_2$ consist of the following 6 triples:
$$ ( c,c,c), (c,c,d),(d,c,c), (d,d,d), (d, d, c), (c, d,d)
 .$$ This
homotopy data essentially differs from $(\alpha_1, S_1)$     by the
choice of $\tau$. The homotopy data $(\alpha_2,
S_2)$ is admissible. As a shift involution $\nu_2$ in $\alpha_2$, we
take the permutation of $c$ and $d$. Nanophrases over $\alpha_2$
considered up to $\nu_2$-permutations and $\nu_2$-shifts of words are
called
  {\it
  quasi-links}.

Connections  to knot theory go as follows. Consider the projection $
\alpha_*\to \alpha_2 $ sending $a_+, b_-$ to $c$ and $a_-,
  b_+$ to $d$. This projection transforms $S_*\subset (\alpha_*)^3$ into $S_2$, commutes with   $\tau$
 and  with
  the shift.
Applying this projection, we can transform any nanophrase over
$\alpha_*$
  into a nanophrase over $\alpha_2$.   This induces a  mapping  ${\PPP}  (\alpha_*,S_* )\to {\PPP}  (\alpha_2,S_2 )$.
  Quotienting by permutations
  and shifts of words (and eventually by inversions of words) we obtain
  projections
  ${\PPP}  (\alpha_*,S_*,
  \nu)\to {\PPP}  (\alpha_2,S_2, \nu_2)$ and
  ${\PPP}_u  (\alpha_*,S_*,
  \nu)\to {\PPP}_u  (\alpha_2,S_2,  \nu_2)$.

A   link
  diagram  $D$    yields a nanophrase over $\alpha_*$ whose projection to
  $\alpha_2$ is
  a nanophrase   over $\alpha_2$ denoted $q(D)$. If $D$ is pointed and ordered, then $q(D)$ is well defined, otherwise
$q(D)$ is defined only up
   to
  $ \nu_2$-permutations and $ \nu_2$-shifts of words. The class of $q(D)$ in ${\PPP}(\alpha_2,S_2,  \nu_2)$ is an invariant of stable
   equivalence of $D$. We call $q(D)$   the {\it underlying quasi-link} of
  $D$. Further projecting to
   ${\PPP}_u (\alpha_2,S_2,  \nu_2)$ we obtain an invariant
   independent of the orientation of $D$.

 Quasi-links will be further discussed in Sect.\ 9.

  \subsection{Remarks}\label{verscurvrelllkiisword}
1. As   explained above,    there are three natural projections from
the set of link diagrams on
  surfaces
  to simpler objects.
  They    map  a link
diagram    to the underlying family of curves,   the underlying
pseudo-link and the underlying quasi-link. The length of the resulting nanophrases   is 
  equal to the number of link components.

 2.
The homotopy data $(\alpha_0, S_0)$ and $(\alpha_1, S_1)$ are
closely related. Consider the bijection from $\alpha_0=\{a,b\}$ to
$\alpha_1=\{1,-1\}$ sending $a$ to $1$ and $b$ to $-1$. This
bijection commutes with the involution $\tau$ in $\alpha_0,
\alpha_1$ and transforms $S_0\subset (\alpha_0)^3$ into a subset of
$S_1 \subset (\alpha_1)^3$. In this way any nanophrase over
$\alpha_0$ determines a  pseudo-link and homotopic nanophrases yield
homotopic pseudo-links. Thus, the homotopy theory of pseudo-links is
a quotient of the homotopy theory of curves. However, there is no way to recover the
pseudo-link $p(D)$ underlying a link diagram $D$  from the system of
curves underlying $D$.
 Note also that the shift involutions in $\alpha_0$ and
$\alpha_1$  do not match: the first one permutes $a$ and $b$ while
the second one is the identity.

  \section{The bracket polynomial}

 The aim  of this section is to construct a polynomial invariant of pseudo-links
  whose value on the underlying pseudo-link of a link diagram is equal to  the Jones polynomial of the link.
   We begin by recalling  Kauffman's   bracket
  polynomial.

   \subsection{Bracket polynomial of links}\label{brack}
L. Kauffman \cite{kau1} associated with   every
  non-empty   link diagram $D$ on a    surface  a 1-variable Laurent polynomial $\langle D\rangle$ called the  {\it
bracket polynomial}  of $D$. This   polynomial is defined by
expanding  each crossing of $D$ as a linear combination of two
uncrossings with coefficients $t$ and $t^{-1}$, see Figure 3. This expands $D$ as
a linear combination of diagrams with no crossings. Each
$d$-component diagram with no crossings is then replaced with
$-(t^2+t^{-2})^{d-1}$.  The bracket polynomial   depends neither on
the orientation of $D$, nor on an order of its components, nor on a
choice of base points. The bracket polynomial is invariant under the
second and the third Reidemeister moves and is multiplied by $t^{\pm
3}$ under the first Reidemeister move.

\begin{figure}
\centerline{\includegraphics[width=9cm]{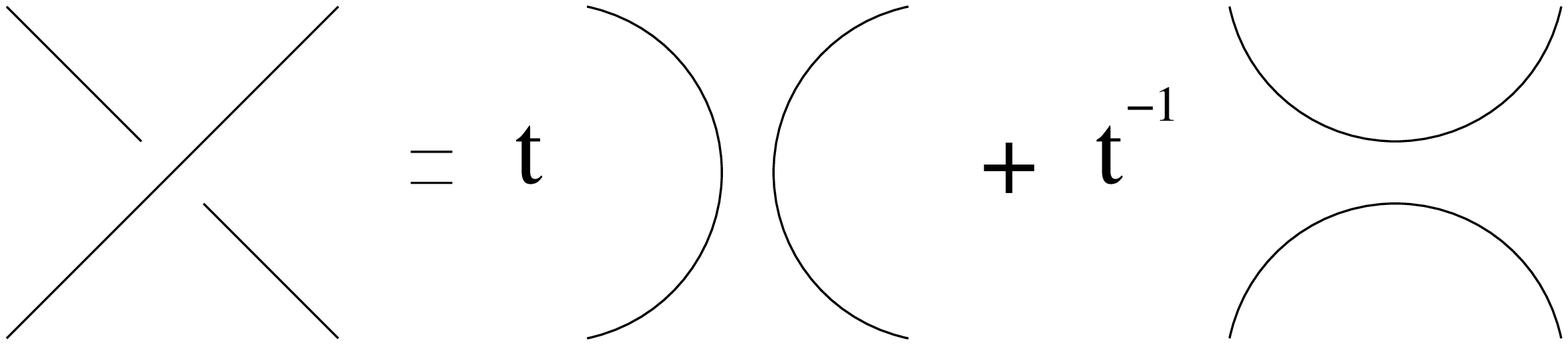}}
\caption{Kauffman's expansion}\label{figure3}
\end{figure}

We  can translate the bracket polynomial to the language of
nanophrases over $\alpha_*$. A nanophrase $P$ over $\alpha_*$ gives
rise to a pointed ordered   link diagram on a  surface. Let $\langle
P\rangle\in \ZZ[t,t^{-1}]$ be the bracket polynomial of this
diagram.  This polynomial is invariant under $\nu$-shifts,
$\nu$-inversions, and $\nu$-permutations on words in $P$ since they
are translated to diagrams as change of base points, orientation
reversal, and change of order of components. The polynomial $\langle
P\rangle$ is invariant under the second and third $S_*$-homotopy
moves on $P$ since they are translated to diagrams as the second and
 third Reidemeister moves. Under the move deleting $AA$ from a
word of $P$, the bracket polynomial is multiplied by $-t^{-3
\varepsilon (A)}$ where $\varepsilon (A)=+1$ if $\vert A\vert \in
\{a_+, b_+\}$ and $\varepsilon (A)=-1$ if $\vert A\vert \in \{a_-,
b_-\}$. When an empty word $\emptyset$ is deleted from a nanophrase
$P$ of length $\geq 2$ the polynomial $\langle P\rangle$ is divided
by $-(t^2+t^{-2})$. Clearly,  $\langle (\emptyset) \rangle=1$.

To    translate the Kauffman crossing expansion to this setting, we
need the following notation.
 Given a  phrase $P$   in an
$\alpha_*$-alphabet $\A$ and a word $w$ in the alphabet $\A$, denote
by $P_w$ the same phrase $P$ in the $\alpha_*$-alphabet $\A_w$
defined in Sect.\ \ref{opernanaopyclomot}. Denote by $w^-$
 the word in the alphabet $\A$ obtained by writing the
letters of $w$ in the opposite order.

The Kauffman crossing expansion applied to  a self-crossing of a
link component and to a crossing of two different components implies
the following two recursive relations for the bracket polynomial of
nanophrases:
$$\langle (P_1 \vert AwAz\vert P_2 ) \rangle=
t^{\varepsilon (A)} \langle (P_1 \vert  w\vert z\vert P_2 )\rangle +
t^{-\varepsilon (A)} \langle (P_1 \vert w^-
 z\vert P_2 )_w\rangle,$$
 $$\langle (P_1 \vert Aw \vert Az\vert P_2 )
\rangle= t^{\varepsilon (A)} \langle (P_1 \vert  w z\vert P_2 )
\rangle +  t^{-\varepsilon (A)} \langle (P_1 \vert w^-
 z\vert P_2 )_w\rangle.$$
 Here  $w$ and $z$ are words in an $\alpha_*$-alphabet $\A$, $A$ is a letter in $ \A$, and
  $P_1,P_2$ are finite sequences
of words  in   $\A$ such  that  every letter of $\A$ appears in the
phrase $(P_1 \vert Aw   Az\vert P_2 )$ twice. In the first formula
$w$ and $z$ are parts of a word $AwAz$  while in the second formula
$Aw$ and $Az$ are two consecutive words. These formulas  and the
properties of $\langle P\rangle$ listed above allow us to compute
this polynomial recursively.

\subsection{Bracket polynomial of pseudo-links}\label{quasibrack}
 In the  recursive formulas above, the right hand side depends
on $\varepsilon (A)$ rather than on $A$. This   implies (by
induction on the number of letters in a nanophrase) that the bracket
  of a nanophrase  over $\alpha_*$ depends only on the
underlying nanophrase  over $\alpha_1$.   Any pseudo-link $p$
determines
  a Laurent
polynomial $\langle p\rangle\in \ZZ[t^{\pm 1}]$ by $\langle
p\rangle=\langle \tilde p\rangle$ where $\tilde p$ is any nanophrase
over $\alpha_*$ whose projection to $\alpha_1$ equals $p$. The
polynomial $\langle p\rangle$   is invariant under shifts,
inversions, and permutations of words in $p$. It is preserved under
the second and third $S_1$-homotopy moves and  is multiplied by
$-t^{-3 \vert A\vert}$ under the move deleting $AA$ from a word of
$p$ (where $\vert A\vert \in \alpha_1=\{1,-1\}$). To compute
$\langle p\rangle$ one can use the   recursive relations   above
with $\varepsilon (A)$ replaced everywhere by $\vert A\vert \in
\{1,-1\}$.

As an illustration, we compute the bracket   for the nanoword
$ABCABC$ over $\alpha_1$ where $\vert A\vert=\vert B\vert=\vert
C\vert  =1$. We have
$$\langle ABCABC \rangle =t\langle (B C \vert BC)\rangle +t^{-1}
\langle ( C_\tau B_\tau B_\tau C_\tau )\rangle$$
$$=t\,(t\langle (CC)\rangle+ t^{-1} \langle (C_\tau C_\tau )\rangle)
+t^{-1} (-t^{3\vert B_\tau \vert }) \langle ( C_\tau C_\tau
)\rangle$$
$$=t^2 (-t^{3\vert C\vert  }) -t^{3\vert  C_\tau \vert }
+t^{-1} (-t^{3\vert B_\tau \vert })(-t^{3\vert C_\tau \vert })
=-t^5-t^{-3}+t^{-7}$$ where $\vert B_\tau \vert =\tau (\vert
B\vert)=  -1$ and  $\vert C_\tau \vert =\tau (\vert C\vert)=  -1$.
This is compatible with the usual formula for  the bracket
  of the standard diagram of a  trefoil.

\subsection{The Jones polynomial of pseudo-links}\label{relllkiisword} For a  pseudo-link $p$ we
 define the
{\it writhe} $\vert  p \vert =\sum_A \vert A\vert\in \ZZ$ where $A$
runs over all letters occurring in $p$. The polynomial $J(p)= (-t)^{-3
\vert p \vert }\langle p\rangle$ is   invariant under all
$S_1$-homotopy moves on $p$. For a pseudo-link $p$ arising from a
  link in
 $    S^3$,    the polynomial  $\langle p\rangle$ is equal to the Jones polynomial of this link   up
to a re-parametrization.

\subsection{Polynomials of phrases}\label{phraseJonbrackword}
An $\alpha_1$-alphabet  is nothing but  a bipartitioned set, that is a set  decomposed as a disjoint union of two    subsets (the preimages of $\pm 1\in \alpha_1$).  Any   phrase  $P$  in an $\alpha_1$-alphabet $\A$ gives rise to a  polynomial 
$\langle P \rangle \in \ZZ[t^{\pm 1}]$ as follows. It is explained in \cite{tu2} that a word $w$  in any alphabet    determines in a canonical way a nanoword  $w^d$ over this alphabet. The same procedure applies to phrases and derives from $P$ a nanophrase $P^d$ over $\A$. (Each letter 
$A\in \A$ occurring $m_A$ times in $P$ gives rise to $m_A(m_A+1)/2$ distinct  letters each occurring in $P^d$ twice.)
Composing with   projection $  \A\to \alpha_1$ we obtain  from $P^d$ a pseudo-link  $(P^d)_1$. Set $\langle P \rangle=\langle (P^d)_1 \rangle$. Similarly, we   define the Jones polynomial of $P$ by
$J(P)=J((P^d)_1)= (-t)^{-3
\vert P \vert } \langle P \rangle$ with   $\vert P \vert=\sum_{A\in \A} \vert A\vert  m_A(m_A+1)/2$ where $\vert A\vert=\pm 1$ is the image of $A$ in $\alpha_1$ and $m_A$ is the number of entries of  $A$ in $P$. The  polynomials $\langle P \rangle, J(P)$ are  interesting invariants of phrases in bipartitioned alphabets.
One  natural question is to characterize the polynomials that arise from  phrases in this way.

 \section{Keis}

Keis were introduced by M. Takasaki  in 1942 as   abstractions of symmetries,  see
\cite{ka}  for  a survey of keis and related objects (quandles, racks, etc.). 
  Here we recall from \cite{tu2} the concept of an $\alpha$-kei
where $\alpha$ is a set with involution $\tau$. This   will be instrumental in the next section 
where we discuss keis of nanophrases.

\subsection{$\alpha$-keis }\label{alphakeisword}
  Consider  a   set
$X$ and suppose that each $a\in \alpha$ gives rise to a bijection
$x\mapsto ax: X\to X$ and to a binary operation $ (x,y)\mapsto
x\ast_a y  $  on $X$. The set $X$   is an {\it $\alpha$-kei} and the
mappings  $ x\mapsto ax ,  (x,y)\mapsto x\ast_a y  $   are {\it kei
operations} if the following axioms are satisfied:

(i) $ax\ast_a x=x$ for all $a\in \alpha, x\in X$;

(ii) $a (x\ast_a y)= ax \ast_a ay$ for all $a\in \alpha, x, y\in X$;

(iii) $ (x\ast_a y)\ast_a z= (x \ast_a az) \ast_a (y\ast_a z)$ for
all $a\in \alpha, x, y, z\in X$;

(iv)      $a\tau (a) x=x$ for all  $x \in X, a\in \alpha$ and

(v)
  $(x\ast_a y) \ast_{\tau(a)} ay=x$  for all  $x, y\in X, a\in \alpha$.

A morphism  of $\alpha$-keis $X\to X'$ is a  set-theoretic mapping  commuting with the kei operations in $X, X'$. 
Given an $\alpha$-kei  $X $, we define an $\alpha$-kei $\overline X$
to be the same  set $X$ with new kei operations $ax: = \tau (a) x$,
$x\ast_a y : =x \ast_{\tau  (a)} y$ for $x,y\in X, a\in \alpha$.
Clearly, $\overline  {\overline X}=X$.

In analogy with group theory, one can define presentations of
$\alpha$-keis by generators and relations.  A presentation of an
$\alpha$-kei  $X$ by generators and relations yields a presentation
of $\overline X$ by generators and relations by replacing every
letter $a\in \alpha$ appearing in the relations by $\tau(a)$.

In the simplest case where $\alpha=\{a\}$ is a 1-element set, an $\alpha$-kei is a set $X$ with involution $x\mapsto \tilde x=ax$ and a binary operation 
$(x,y)\mapsto x\ast  y=  x\ast_a y$ such that 
$\tilde x \ast  x=x$;
 $\widetilde {x\ast  y}= \tilde x \ast \tilde y$;
  $ (x\ast  y)\ast  z= (x \ast \tilde z) \ast  (y\ast z)$, and $(x\ast  y) \ast \tilde y=x$ for
all $  x, y, z\in X$. When the involution $x\mapsto \tilde x$ is the identity, these axioms are equivalent to those of  a kei, see \cite{ka}.

The  $\alpha$-keis generalize  quandles: there is a canonical 
bijection  between quandles  and  $\alpha$-keis $X$ such that  $\alpha=\{a,b\}$ is a
2-element set with involution permuting $a,b$   and  $ax=bx=x$
for all $x\in X$, cf.\  \cite{tu2}, Lemma 14.7.1. 

\subsection{Core $\alpha$-keis}\label{hoexamplsword} The following construction  of $\alpha$-keis provides a vast set of examples. By a {\it $\tau$-compatible action} of $\alpha$ on a group $G$  we   mean   a set of  group automorphisms $\{G\to G, g\mapsto ag\}_{a\in \alpha}$ such that $a \tau (a) g=g$ for all $a\in \alpha, g\in G$. It is easy to check that such an action   together with   kei operations $g\ast_a h= h (\tau (a) g)^{-1} h$ make $G$ into an $\alpha$-kei. It is called the {\it core} of $G$ and  denoted ${\rm {core}} (G)$. For $\alpha$ consisting of one  element that acts on $G$ as the identity, this construction is due to  D.\ Joyce (cf.\  \cite{fr}, p.\ 349).

The construction of the core has a natural adjoint  associating with an arbitrary $\alpha$-kei $X$ a group $\Gamma_X $ 
with generators $\{[x]\}_{x\in X}$ and relations  $[a(x\ast_b y)]=[a y] [a  \tau (b) x]^{-1} [ay]$ for all $a,b\in \alpha, x,y \in X$.  We endow $\Gamma_X$ with the $\tau$-compatible action of  $\alpha$   defined on the generators by $a[x]=[ax]$ for $a\in \alpha, x\in X$.   Given a group $G$ with $\tau$-compatible action of $\alpha$  and a kei morphism
$f:X \to  {\rm {core}} (G)$, there is a unique group homomorphism $ \Gamma_X\to G$ whose composition   with the inclusion $X\hookrightarrow  \Gamma_X, x\mapsto [x]$ is equal to $f$. This universal  property characterizes $\Gamma_X$  up to isomorphism.

A presentation of $\Gamma_X$ by generators and relations can be read from an arbitrary presentation $[S:R]$ of $X$ by generators and relations (cf.\  \cite{fr}, Lemma 4.3). Namely,  $\Gamma_X$ is generated by the symbols $\{as\}_{a\in \alpha, s\in S}$ subject to the relations obtained from $R$ by replacing all terms 
of type  $a(x\ast_b y) $ by $(a y) (a  \tau (b) x)^{-1}   (ay)$.

   \section{Keis of nanophrases}

A    well known construction due to S. Matveev and D. Joyce associates quandles with link diagrams.  Since link diagrams are  nanophrases over  $\alpha_*$,  one may attempt to generalize this construction to   nanophrases
over
an arbitrary  alphabet $\alpha$   with involution $\tau$. We do it  here starting with certain additional data.

\subsection{Keis  and homotopy}\label{homotalphakeisword}  Fix
an equivalence relation $\sim$ on $\alpha$ such  that   $a\sim b \Rightarrow \tau(a)\sim
\tau(b)$ for $a,b \in \alpha$. Let   $\underline \alpha=\alpha/\sim$
with involution  $\underline \tau$ induced by $\tau$.
 For   $a\in \alpha$, denote its projection to $\underline \alpha$ by
  $\underline a$.

Fix
  a set (possibly empty) $\beta \subset \alpha $    such that
  $\tau (\beta)=\beta$. We associate with any nanophrase $P=(\A,  (w_1,...,w_k))$ over
$\alpha$ an $\underline \alpha$-kei $   \kappa_{\beta } (P)$ as follows. Let $n_r$ be the
length of the word $w_r$ for $r=1,...,k$.  Each letter  $A\in \A$
appears in $P$ twice, say,  first time at the $i_1$-th position in
$w_{r_1}$ and second time at the $i_2$-th position  in $w_{r_2}$
where $1\leq i_1\leq n_{r_1}, 1\leq i_2\leq n_{r_2}$, $r_1\leq r_2$,
and $r_1=r_2 \Rightarrow i_1 < i_2$. The $\underline \alpha$-kei $\kappa_{\beta }
(P) $   is generated by the symbols
$\{x^r_s\}$ where $1\leq r \leq k$ and $0\leq s \leq n_r$.  Each $A\in \A$
gives rise to two defining relation: if $a=\vert A\vert \in \beta$, then
$$x_{i_1}^{r_1}=\underline a \, x_{i_1-1}^{r_1}, \,\,\,\,\,x_{i_2}^{r_2}=x_{i_2-1}^{r_2} \ast_{\underline  a} x_{i_1-1}^{r_1},$$
and if  $a=\vert A\vert \in \alpha - \beta $, then
$$x_{i_1 }^{r_1}=x_{i_1-1}^{r_1}\ast_{\underline  a} x_{i_2-1}^{r_2} ,  \,\,\,\,\,x_{i_2}^{r_2}
=\underline a  \,x_{i_2-1}^{r_2} .$$
 These generators and relations define the $\underline \alpha$-kei $\kappa_{\beta }
(P)$. It  has two sets of distinguished elements
   $x^1_0, x^2_0,..., x^k_0$  (the  inputs)  and  $x^1_{n_1}, x^2_{n_2},..., x^k_{n_k}$ (the    outputs).
 Adding  the relations  $x^r_0=x^r_{n_r}$  for
$r=1,...,k$ we obtain a quotient  $\underline \alpha$-kei $  \hat{\kappa}_{\beta } (P) $.

  Note the obvious   $\underline
\alpha$-kei isomorphism  $\kappa_{\beta } (P^-)\approx
\overline {\kappa_{\alpha-\beta } (P)}$ where $P^-$ is $P$ read from
right to left. This isomorphism transforms the $r$-th input (resp.\
output) into the $(n+1-r)$-th output (resp.\ input).
 Clearly,  $\hat
{\kappa}_{\beta } (P^-)\approx \overline {\hat {\kappa}_{\alpha-\beta} (P)}$.

For the next definition, it is convenient to set $\beta_0=\beta$ and $\beta_1=\alpha-\beta$. Let $S=S(\beta, \sim)\subset \alpha^3$
consist  of all triples $(a,b,c)\in \alpha^3$ such that

-   $a\sim b\sim c$ and $a,b,c\in \beta_i$ for some $i\in \{0,1\}$;

- or $a\sim b\sim \tau (c)$ and $a,b\in \beta_i, c\in \beta_{1-i}$
for some $i\in \{0,1\}$;

- or   $\tau (a)\sim b\sim c$ and $b,c\in \beta_i, a\in \beta_{1-i}$
for some $i\in \{0,1\}$.

The set $S$ contains the diagonal of $\alpha^3$ and therefore the
homotopy data $(\alpha, S)$ is admissible.

\begin{theor}\label{erddcppmolbrad}  For any   nanophrase $P $ over
$\alpha$, the $\underline \alpha$-kei $    \kappa_{\beta } (P)$ is invariant under $S$-homotopy moves.
\end{theor}

The proof   goes by repeating the  proof of Lemma 15.1.1 in
\cite{tu2}.

  The next theorem shows that for an
appropriate choice of the shift involution $\nu$, the kei $\kappa_{\beta } (P)$   is also
preserved under $\nu$-permutations on $P$ and its quotient $\hat
{\kappa}_{\beta } (P)$ is preserved under $\nu$-permutations and $\nu$-shifts.

\begin{theor}\label{2223erddcppmolbrad} Let $\nu:\alpha\to \alpha$ be an involution  such that
$\nu(\beta)=\alpha- {\beta }  $ and $a\sim \nu (a)$ for all $a\in \alpha$.
For a    nanophrase $P $ over $\alpha$, the $\underline \alpha$-kei
$   \kappa_{\beta } (P)$ is invariant under
$\nu$-permutations on the words of $P$. The quotient $\underline
\alpha$-kei $\hat {\kappa}_{\beta } (P) $  is invariant under  $\nu$-permutations
and $\nu$-shifts  on the words of $P$.
\end{theor}

We leave the proof to the reader as an exercise.

\subsection{Examples}\label{fronslalphakeisword} 1. Consider the
alphabet $\alpha_*=\{a_+, a_-, b_+, b_-\}$ with involution
$\tau(a_\pm)=b_\mp$,   shift involution $\nu(a_\pm)=b_\pm$, and
  distinguished subset 
$\beta =\{a_+,b_-\}$. Provide $\alpha_*$
with   equivalence relation $a_+\sim b_+, a_-\sim b_-$. This data
satisfies all the conditions of Theorems \ref{erddcppmolbrad} and
\ref{2223erddcppmolbrad}. Clearly,  $S(\beta ,
\sim)=S_*\subset (\alpha_*)^3$  is the set defined in
Sect.\ \ref{triplerss}. This yields for   any nanophrase $P$   over
$\alpha_*$    an $\underline \alpha_*$-kei $ 
 {\kappa}_{\beta } (P) $ invariant under $S_*$-homotopy. The quotient $\underline \alpha_*$-kei $\hat{\kappa}_{\beta } (P)$ is also  invariant under  $\nu$-shifts and $\nu$-permutations. The set $\underline
\alpha_*=\{+,-\}$ consists   of 2 elements permuted by $\underline \tau$.
For the nanophrase $P$ associated with a link diagram on a    surface,
the $\underline \alpha_*$-kei $\hat{\kappa}_{\beta } (P)$ is invariant
under stable equivalence and independent of the choice of the order
and the base points of the link components. The quotient of $\hat{\kappa}_{\beta } (P)$ by $ax=x$ for all $a\in \alpha_*=\{+,-\}, x\in \hat{\kappa}_{\beta } (P)$   with binary  operation $\ast_+$ is the standard link quandle (see \cite{fr}, \cite{ka}, \cite{kau}).

2. Consider the alphabet $\alpha_0=\{a, b\}$ with involution $\tau(a
)=b $  and   distinguished subset
  $\beta_0=\alpha $. As the    equivalence relation $\sim$  in
$\alpha_0$ we take the equality $=$. This  data satisfies  the
conditions of Theorem  \ref{erddcppmolbrad} where $S(\beta_0 ,
\sim)=S_0\subset (\alpha_0)^3$  is the diagonal. For any nanophrase $P$   over $\alpha_0$, we obtain    an $\underline
\alpha_0$-kei $   {\kappa}_{\beta_0 } (P) $
invariant under $S_0$-homotopy.   This
example is contained in  the previous one: the mapping $\alpha_0\to \alpha_\ast$ defined by $a\mapsto a_+, b\mapsto b_-$ transforms $P$ into a nanophrase $ P_*$
  over $\alpha_*$  and  $ {\kappa}_{\beta_0 } (P)={\kappa}_{\beta } (  P_*)$. 
   
3. The homotopy data $(\alpha_1, S_1)$  from  Sect.\ \ref{erdt} cannot be obtained by the methods of Sect.\
\ref{homotalphakeisword} and does not lead to   
  keis. Pseudo-links have no keis.

  4. Consider the alphabet $\alpha_2=\{c, d\}$ with  trivial involution $\tau=\id $,
  shift involution $\nu $ permuting $c$ and $d$,
   and  distinguished subset 
  $\beta =\{c\}$. Provide  $\alpha_2$ with trivial  equivalence
relation $\sim$ (all elements are equivalent). This data satisfies
  the conditions of Theorems \ref{erddcppmolbrad} and
\ref{2223erddcppmolbrad}. Clearly, $S(\beta ,
\sim)=S_2\subset (\alpha_2)^3$ is the set defined in
Sect.\ \ref{quasinanaopyclomot}. This yields for any quasi-link
 $P$ (= a nanophrase   over $\alpha_2$) an $\underline \alpha_2$-kei
$ {\kappa}_{\beta } (P) $ invariant under
$S_2$-homotopy. The quotient  $\underline \alpha_2$-kei
$\hat{\kappa}_{\beta } (P)$ is  also  invariant   under  $\nu$-shifts and $\nu$-permutations. The set
$\underline \alpha_2$ consists here   of 1 element. A study of this
kei should lead to interesting homotopy invariants of quasi-links. 

When $P$ is obtained from 
an oriented link   $L\subset S^3$ by taking the associated nanophrase over $\alpha_*$ and projecting to $\alpha_2$, the    group 
$\Gamma=\Gamma_{\hat{\kappa}_{\beta } (P)}$ (defined in Sect.\ \ref{hoexamplsword}) is closely related to the fundamental group of the 2-fold branched cover $M$ of $S^3$ with branching set  $L$. Namely, the group $\pi_1(M) \ast \ZZ$ is the quotient of 
$\Gamma $ by the relations $ag=g$ for the unique $a\in \underline \alpha_2$ and all $g\in \Gamma $. This is obtained by comparing  the presentation of ${\hat{\kappa}_\beta (P)}$ as above with  the Wada    presentation of  $\pi_1(M) \ast \ZZ$, both computed  
  from a diagram of $L$. This observation easily extends to the nanophrase derived from (a  diagram of) a link $L\subset   \Sigma \times [0,1]$ where $\Sigma$ is a surface. Here one should use  the 
  2-fold branched cover   of $ \Sigma \times [0,1]/ \Sigma \times 1$ with branching set  $L$, cf.\ the argument  in \cite{kk}, Prop.\ 5.1 and the proof of Wada's theorem in \cite{pr}, p.\ 287.

5. Let $\alpha_*, \tau, \nu, \beta ,\sim$ be as in Example 1.  Pick a set $\gamma$ and consider the
alphabet $\alpha_\gamma=\alpha_*\times \gamma$  with involution
$\tau\times \id$,   shift involution $\nu_\gamma=\nu \times \id$,  and distinguished subset
  $ \beta \times \gamma$. Provide $\alpha_\gamma$ with equivalence relation $\sim_\gamma$ as follows: two pairs  $(x, c),(y,d)$ with $x,y\in \alpha_*, c,d \in \gamma$ are equivalent if    $x\sim y$ and $c=d$.
This data
satisfies   the conditions of Theorems \ref{erddcppmolbrad} and
\ref{2223erddcppmolbrad}. The set  $S_\gamma=S(\beta  \times \gamma,
\sim_\gamma)$ is the product of $S_*\subset (\alpha_*)^3$ and the diagonal of $\gamma^3$. This yields for   any nanophrase $P$   over
$\alpha_\gamma$    an $\underline {\alpha_\gamma}$-kei $ 
 {\kappa}_{\beta \times \gamma} (P) $ invariant under $S_\gamma$-homotopy where $\underline {\alpha_\gamma}=\alpha_\gamma/\sim_\gamma= \{+,-\}\times \gamma$. The quotient $\underline {\alpha_\gamma}$-kei $\hat{\kappa}_{\beta \times \gamma} (P)$ is also  invariant under  $\nu_\gamma$-shifts and $\nu_\gamma$-permutations.

Nanophrases over $\alpha_\gamma$ have a simple geometric interpretation.   Let us call an (ordered pointed) link diagram  on an (oriented) surface  {\it $\gamma$-colored} if all
its crossings are endowed with elements of $\gamma$ (the colors).
Homeomorphisms of $\gamma$-colored link diagrams should preserve the colors
of the crossings. Stable equivalence of $\gamma$-colored link diagrams is
defined as in the non-colored case with the following restrictions on the Reidemeister
moves:
  the second move is allowed only when it involves two
crossings of the same  color,  the third   move is allowed only when
it involves three crossings of the same  color  which is kept under
the  move. The crossings not involved in the moves   keep their
color. Denote $\mathcal L_\gamma$ the set of stable equivalence classes of
$\gamma$-colored ordered pointed link diagrams. The same arguments as in Theorem \ref{rtlinkslbrad}  show that $\mathcal L_\gamma= \PPP(\alpha_ \gamma, {S_\gamma})$. This equality
  implies a 
similar equality for non-ordered   non-pointed link diagrams. The construction  above associates  with every $\gamma$-colored link diagram an   $\underline {\alpha_\gamma}$-kei invariant under stable equivalence.

\subsection{Remark}\label{alph300029word} Further   invariants of  nanophrases   can be
derived from their keis  by abelianization   \cite{tu2}. Another
interesting possibility   is to define cohomology of $\alpha$-keis and to derive
homotopy invariants of  nanophrases from cocycles and
state sums.

                     \end{document}